\documentclass[11pt]{amsart}
\usepackage{palatino}
\usepackage{amsfonts}
\usepackage{amssymb}
\usepackage{amscd}

\usepackage{amsmath}
\usepackage{mathrsfs}

\newtheorem{theorem}{Theorem}[section]
\newtheorem{proposition}[theorem]{Proposition}
\newtheorem{corollary}[theorem]{Corollary}
\newtheorem{lemma}[theorem]{Lemma}
\theoremstyle{definition}

\def\OO{\mathcal{O}}
\def\E{\mathcal{E}}
\def\A{\mathcal{A}}
\def\U{\mathcal{U}}
\def\L{\mathcal{L}}
\def\Q{\mathcal{Q}}
\newcommand{\HH}{\mathcal{H}}

\newcommand{\PP}{\mathbb{P}}

\newcommand{\rk}{\operatorname{rk}}

\begin{document}
\title[Semipositivity and complete subvarieties]
{Semipositivity and complete subvarieties of projectivised strata of holomorphic differentials}

\author[S. Mullane]{Scott Mullane}
\address{School of Mathematics and Statistics, University of Melbourne, Australia}
\email{{\tt mullanes@unimelb.edu.au}}

\begin{abstract}
We show that if the projectivised stratum $\PP\HH_g(\mu)$ associated to a partition $\mu=(m_1,\ldots,m_n)$ of $2g-2$ contains a complete curve, then
\[
\kappa_\mu:=\sum_{i=1}^n\frac{m_i(m_i+2)}{m_i+1}>12.
\]
The weak inequality follows from semipositivity applied to a quotient of the Hodge bundle, while equality is excluded using the second Fujita decomposition and the geometry of fibred surfaces with maximal unitary rank.  Consequently, no projectivised holomorphic stratum in genus $g\leq5$ contains a positive-dimensional complete subvariety. In genus $g=6$ the same conclusion holds for
\[
(10),\ (9,1),\ (8,2),\ (8,1,1),\ (7,3),\ (6,4),\ (5,5).
\]
\end{abstract}

\maketitle


\section{Introduction}

Let $g\geq2$ and let $\mu=(m_1,\ldots,m_n)$ be a partition of $2g-2$ with all $m_i>0$. The stratum $\HH_g(\mu)$ parametrises the data $(C,\omega,p_1,\ldots,p_n)$ consisting of a smooth curve of genus $g$ and a nonzero holomorphic differential whose zero orders are prescribed by $\mu$ and occur at the marked points.  Its projectivisation $\PP\HH_g(\mu)$ parametrises the differential only up to scale.  Throughout, the zeros are labelled and ordered.  Forgetting the ordering gives a finite map and does not affect the existence of complete subvarieties. 

The unprojectivised strata contain no complete curves by~\cite{Gendron}, see also Chen's algebro-geometric proof in~\cite{ChenComplete}.  Projectivisation removes the scaling function used in those arguments and changes the problem substantially.  The existence of complete curves in projectivised holomorphic strata remains open in general, even in low genus.  See~\cite[Section~6.2]{ChenYu}.

For strictly meromorphic differentials, Chen proved that projectivised strictly meromorphic strata contain no positive-dimensional complete subvarieties~\cite{ChenAffineGeometry,ChenAffine}, and Chen and Tahar subsequently gave a flat-geometric proof of this result~\cite{ChenTahar}. Their result does not apply to the holomorphic strata considered here.

For a signature $\mu$ define
\begin{equation}\label{eq:kappamu}
 \kappa_\mu:=\sum_{i=1}^n\frac{m_i(m_i+2)}{m_i+1}
 =2g-2+n-\sum_{i=1}^n\frac{1}{m_i+1}.
\end{equation}
Our main result gives a strict numerical obstruction at $\kappa_\mu=12$.

\begin{theorem}
Let $\mu$ be a signature of holomorphic differentials.  If $\PP\HH_g(\mu)$ contains a complete curve, then
\[
\kappa_\mu>12.
\]
Consequently, if $\kappa_\mu\leq 12$, then $\PP\HH_g(\mu)$ contains no positive-dimensional complete subvariety.
\end{theorem}

The numerical obstruction immediately gives the following.

\begin{corollary}\label{cor:strata}
Every projectivised stratum of holomorphic differentials in genus $g\leq 5$ contains no positive-dimensional complete subvariety. If $g=6$  then $\PP\HH_6(\mu)$ contains no positive-dimensional complete subvariety for $\mu$ equal to
\[
(10),\ (9,1),\ (8,2),\ (8,1,1),\ (7,3),\ (6,4),\ (5,5).
\]
\end{corollary}

The argument has two parts.  Let $Z$ be an integral complete curve in a projectivised stratum.  After normalising $Z$ and making a finite base change, there is a smooth family of curves $f:X\longrightarrow B$, with Hodge bundle
\[
\E=f_*\omega_f,
\]
and a tautological line subbundle $\L\subset\E$ generated fibrewise by the distinguished differential.  The standard tautological relation on the stratum becomes
\begin{equation}\label{eq:introrelation}
 12\deg\E=\kappa_\mu\deg\L.
\end{equation}
Chen's tautological calculation and the relative Noether formula give~\eqref{eq:introrelation}, see Proposition~\ref{prop:relation}.  The Hodge bundle $\E$ is nef by Fujita's semipositivity theorem~\cite[Theorem~(2.7)]{Fujita}, and hence so is its quotient $\E/\L$.  Since $\deg\L>0$, equation~\eqref{eq:introrelation} gives $\kappa_\mu\geq12$.

If equality holds, the quotient $\E/\L$ is nef of degree zero.  The second Fujita decomposition
\[
\E=\A\oplus\U
\]
has $\A$ ample and $\U$ unitary flat.  The map from $\A$ to $\E/\L$ vanishes, so $\U$ surjects onto a bundle of rank $g-1$.  Nonisotriviality then forces the unitary rank to be $g-1$.  For $g\geq4$ such a fibration cannot be a Kodaira fibration by~\cite[Theorem 3.1]{Stoppino}, contradicting the fact that $f$ is smooth and nonisotrivial.  Finally, observe that numerically that equality occurs for in \eqref{eq:kappamu} for exactly one signature, that is, $\mu=(1^8)$ in genus $g=5$.

Affineness implies the absence of positive-dimensional complete subvarieties.  Some of the strata covered by Corollary~\ref{cor:strata} were already known to be affine, including every hyperelliptic component and the low genus nonvarying strata treated in~\cite{ChenAffine}.  However, the new result covers many strata not covered previously. In particular, it includes the principal stratum in genus $g=3,4,5$, and many further components not covered by the known affineness results.

Devkota and Lee show that the rational Chow rings of the multi-scale compactifications of the genus-three strata with signatures $(4)$, $(2,2)$ and $(3,1)$ are generated by boundary-stratum classes~\cite{DevkotaLee}.  This gives another proof that these open strata contain no positive-dimensional complete subvarieties that was already known from affineness.  Their result does not treat the principal stratum $\PP\HH_3(1^4)$.

More precisely, in genus $g=3$ our new case is the principal stratum $\PP\HH_3(1^4)$.  In genus $g=4$ the new cases are
\[
\mu=(4,1^2), (3,1^3), (2^3)^{\rm even}, (2^2,1^2), (2,1^4), (1^6)
\]
In genus $g=5$ it covers every component, including the principal stratum $\PP\HH_5(1^8)$, where $\kappa_\mu=12$.  We use the component labels given by Kontsevich--Zorich~\cite{KZ}.

\medskip
\noindent\textbf{AI disclosure.} Generative AI tools were used at points in this project. The author takes full responsibility for all content.
\smallskip

\noindent\textbf{Acknowledgements.} This project was supported by DECRA Grant DE220100918 from the Australian Research Council.

\section{Families in a projectivised stratum}\label{sec:background}

We begin with the geometry used throughout the paper.  We use the convention that $\PP\HH_g(\mu)$ parametrises one-dimensional subspaces of holomorphic differentials.  Thus its tautological line bundle is the restriction of $\OO_{\PP\E}(-1)$.

We make the standard reduction from such a curve to a family of pointed curves as follows.  Let
\[
 Z\subset\PP\HH_g(\mu)
\]
be an integral complete curve with normalization $B_0\longrightarrow Z$. Choose a full symplectic level structure of level $\ell\geq3$.  The pullback of the stratum to the fine level-$\ell$ moduli space is finite over $\PP\HH_g(\mu)$ and carries the universal pointed curve and the tautological differential line.  Let $B$ be the normalisation of an irreducible component of its fibre product with $B_0$ which dominates $B_0$. Then $B$ is a smooth complete curve, the map $B\longrightarrow B_0$ is finite and surjective, and we obtain a family
\[
 f:X\longrightarrow B
\]
of smooth genus $g$ curves, disjoint sections $\sigma_1,\ldots,\sigma_n$, and a line subbundle
\[
 \L\subset\E:=f_*\omega_f.
\]
Fibrewise, $\L_b$ is the line spanned by a differential with divisor $\sum_i m_i\sigma_i(b)$.  The tautological inclusion is everywhere nonzero and has constant rank one, so $\L$ is a line subbundle and $\E/\L$ is locally free. The induced map from $B$ to the stratum is nonconstant because it is surjective onto $Z$.  A further finite base change multiplies all degrees by the same positive integer, so none of the inequalities below are affected.

Write
\[
 \lambda_B:=\deg\E,\qquad \eta_B:=\deg\L.
\]
The notation is compatible with the tautological classes $\lambda$ and $\eta=c_1(\OO_{\PP\E}(-1))$ on the stratum.

\begin{proposition}\label{prop:relation}
For the family $f:X\longrightarrow B$ constructed above,
\[
 12\lambda_B=\kappa_\mu\eta_B.
\]
\end{proposition}

\begin{proof}
We first recall Chen's tautological calculation \cite[proof of Proposition~2.1]{ChenTautological}.  The universal differential gives an isomorphism
\begin{equation}\label{eq:canonicalfamily}
 \omega_f\cong f^*\L\otimes \OO_X\left(\sum_{i=1}^n m_i\sigma_i\right).
\end{equation}
Restricting~\eqref{eq:canonicalfamily} to $\sigma_i$ and using adjunction gives
\[
 c_1(\L)=(m_i+1)\psi_i
\]
on $B$, where $\psi_i=c_1(\sigma_i^*\omega_f)$ and
$\sigma_i^2=-\deg\psi_i$.  Consequently,
\[
 \sigma_i^2=-\frac{\eta_B}{m_i+1}.
\]
The sections are pairwise disjoint.  Squaring~\eqref{eq:canonicalfamily} therefore gives
\begin{align*}
 K_f^2
 &=2\left(\sum_i m_i\right)\eta_B+
   \sum_i m_i^2\sigma_i^2\\
 &=\sum_i\left(2m_i-\frac{m_i^2}{m_i+1}\right)\eta_B\\
 &=\sum_i\frac{m_i(m_i+2)}{m_i+1}\eta_B
 =\kappa_\mu\eta_B.
\end{align*}
Since $f$ is smooth, the relative Noether formula is $K_f^2=12\deg f_*\omega_f=12\lambda_B$. 
The result follows.
\end{proof}

\begin{lemma}\label{lem:positiveeta}
For the family associated above to an integral complete curve $Z\subset\PP\HH_g(\mu)$, it holds that $\eta_B>0$.
\end{lemma}

\begin{proof}
The Hodge bundle $\E$ is nef by Fujita semipositivity~\cite{Fujita} and hence $\lambda_B\geq0$.  Proposition~\ref{prop:relation} and $\kappa_\mu>0$ imply $\eta_B\geq0$.

Suppose that $\eta_B=0$.  Then $\lambda_B=0$.  A smooth family of curves with Hodge bundle of degree zero is locally trivial, so after a finite unramified base change the family is a product $C\times B\longrightarrow B$.  Hence $B\subset \PP H^0(C,K_C)$. But then as $\L$ is the pullback of $\OO(-1)$, we have $\deg\L<0$ because $\OO(1)$ is ample providing a contradiction to $\eta_B=0$.
\end{proof}

\section{The semipositivity obstruction}

The basic obstruction is as follows.

\begin{theorem}\label{thm:weak}
If $\PP\HH_g(\mu)$ contains a complete curve, then
\[
 \kappa_\mu\geq12.
\]
\end{theorem}

\begin{proof}
We apply the construction of Section~\ref{sec:background} to an integral complete curve $Z\subset\PP\HH_g(\mu)$.  Since the Hodge bundle $\E$ is nef, every quotient of $\E$ is nef~\cite[Proposition~6.1.2(i)]{Lazarsfeld}.  In particular,
\[
 \Q:=\E/\L
\]
is nef, and hence $\deg\Q\geq0$.  Proposition~\ref{prop:relation} gives
\[
 \deg\Q=\lambda_B-\eta_B =\left(\frac{\kappa_\mu}{12}-1\right)\eta_B.
\]
By Lemma~\ref{lem:positiveeta}, $\eta_B>0$ and hence $\kappa_\mu\geq12$.
\end{proof}

\section{The case $\kappa_\mu=12$}

We extend Theorem~\ref{thm:weak} to show the case $\kappa_\mu=12$ cannot occur.  Recall the second Fujita decomposition of the Hodge bundle of a fibration over a curve:
\begin{equation}\label{eq:fujita}
    \E=\A\oplus\U,
    \end{equation}
    where $\A$ is ample and $\U$ is unitary flat.  This is Fujita's second decomposition~\cite{FujitaSecond}, for a complete proof, see~\cite{CataneseDettweilerFull}.  Recall
\[
 u_f:=\rk\U
\]
for the unitary rank of the fibration.

\begin{lemma}\label{lem:equalityrank}
Let $Z\subset\PP\HH_g(\mu)$ be an integral complete curve.  If $\kappa_\mu=12$, then the smooth family constructed in Section~2 has
\[
 u_f=g-1.
\]
\end{lemma}

\begin{proof}
Proposition~\ref{prop:relation} implies $\deg\E=\deg\L$, and hence the nef bundle $\Q=\E/\L$ has degree zero. A nef vector bundle of degree zero on a smooth complete curve is numerically flat. That is, all its Harder--Narasimhan slopes are nonnegative, and their sum is zero, so every slope is zero.  In particular, $\Q$ is semistable of degree zero and $\Q^\vee$ is nef~\cite[Theorem~6.4.15]{Lazarsfeld}.

Consider the quotient map $\E\longrightarrow\Q$ and the decomposition \eqref{eq:fujita}.  There is no nonzero morphism $\A\longrightarrow\Q$. The image of such a morphism would be a nonzero quotient of the ample bundle $\A$, which is an ample bundle of positive degree.  Its saturation in $\Q$ would have at
least the same degree, whereas every subbundle of the semistable degree zero bundle $\Q$ has nonpositive degree.  Thus the quotient map factors through a
surjection
\[
 \U\longrightarrow\Q.
\]
It follows that $u_f\geq\rk\Q=g-1$.  The family is nonisotrivial by Lemma~\ref{lem:positiveeta}.  More precisely, the second Fujita decomposition
and $\kappa_\mu=12$ give
\[
 \deg\A=\deg\E=\lambda_B=\eta_B>0.
\]
Thus $\A$ is nonzero and $u_f\leq g-1$.  Hence $u_f=g-1$.
\end{proof}

\begin{theorem}\label{thm:equality}
Let $g\geq4$.  There is no complete curve in a projectivised holomorphic stratum satisfying $\kappa_\mu=12$.
\end{theorem}

\begin{proof}
Suppose that such a curve exists.  Its associated fibration is smooth and nonisotrivial, and hence is a Kodaira fibration.  By Lemma~\ref{lem:equalityrank} it has unitary rank $u_f=g-1$. However, a relatively minimal nontrivial fibration of genus $g\geq4$ with $u_f=g-1$ is not a Kodaira fibration by~\cite[Theorem 3.1]{Stoppino}. In particular, it satisfies
\[
 K_f^2<12\chi_f.
\]
This contradicts the relative Noether equality $K_f^2=12\chi_f$ for the smooth family $f$.
\end{proof}

Combining Theorems~\ref{thm:weak} and~\ref{thm:equality} gives the main result.

\begin{theorem}\label{thm:strict}
If $\PP\HH_g(\mu)$ contains a complete curve, then $\kappa_\mu>12$.
\end{theorem}

\section{Numerical consequences}

All that is left is to determine the holomorphic signatures to which Theorem~\ref{thm:strict} applies.  Set
\[
 h(m):=\frac{m(m+2)}{m+1}=m+1-\frac{1}{m+1}.
\]
If $a,b\geq1$, then
\begin{equation}\label{eq:splitting}
 h(a)+h(b)-h(a+b)
 =1-\frac{1}{a+1}-\frac{1}{b+1}+\frac{1}{a+b+1}>0.
\end{equation}
Thus splitting a part of a partition strictly increases $\kappa_\mu$.

\begin{proposition}\label{prop:numerics}
Let $\mu$ be a partition of $2g-2$.
\begin{enumerate}
\item If $g\leq4$, then $\kappa_\mu<12$.
\item If $g=5$, then $\kappa_\mu\leq12$, with equality precisely for
$\mu=(1^8)$.
\item If $g=6$, then $\kappa_\mu<12$ precisely for
\[
(10),\ (9,1),\ (8,2),\ (8,1,1),\ (7,3),\ (6,4),\ (5,5).
\]
There is no genus six signature with $\kappa_\mu=12$.
\item If $g\geq7$, then $\kappa_\mu>12$ for every signature $\mu$.
\end{enumerate}
\end{proposition}

\begin{proof}
By~\eqref{eq:splitting}, among partitions of $2g-2$ the principal partition $(1^{2g-2})$ has maximal value
\[
 \kappa_{(1^{2g-2})}=3g-3,
\]
while the minimal partition $(2g-2)$ has value
\[
 \kappa_{(2g-2)}=2g-1-\frac{1}{2g-1}.
\]
Hence parts $(i),(ii)$ and $(iv)$ follow immediately and it remains only to consider $g=6$.  Since $\sum_i m_i=2g-2=10$, formula \eqref{eq:kappamu} becomes
\[
 \kappa_\mu=10+\sum_i\frac{m_i}{m_i+1}.
\]
Set $q(m)=m/(m+1)$ which satisfies
\[
 q(m+1)-q(m)=\frac{1}{(m+1)(m+2)},
\]
and is hence strictly decreasing in $m$.  Consequently, if $a\geq b\geq2$, transferring one unit from $b$ to $a$ strictly decreases $q(a)+q(b)$.  Among partitions of ten with at least four parts, equation~\eqref{eq:splitting} reduces the minimum to the four-part case, where the transfer argument shows that it is attained at $(7,1,1,1)$.  Hence $\kappa_\mu>12$ for every such partition.  Every one- or two-part partition has $\kappa_\mu<12$.

For three parts, the minimum is attained at $(8,1,1)$ and is less than twelve.  Starting from any other three-part partition $(a,b,c)$ with $a\geq b\geq c$, first transfer units from $c$ to $a$ until $c=1$, and then from $b$ to $a$ until $b=2$.  This path never meets $(8,1,1)$ and terminates at $(7,2,1)$.  Hence the minimum among all other three-part partitions is attained at $(7,2,1)$ where $\kappa_\mu>12$. Hence only $(8,1,1)$ occurs.  The inequalities are strict throughout, hence there is no equality case in genus six. Combining this with the preceding bounds shows that the equality $\kappa_\mu=12$ occurs only for $g=5$ and $\mu=(1^8)$.
\end{proof}

\begin{corollary}\label{cor:lowgenus}
For $g\leq5$ and every holomorphic signature $\mu$ of genus $g$, the stratum $\PP\HH_g(\mu)$ contains no complete curve.
\end{corollary}

\begin{proof}
This follows from Proposition~\ref{prop:numerics} and Theorem~\ref{thm:strict}.
\end{proof}

\begin{corollary}\label{cor:g6}
The projectivised genus $g=6$ strata with signatures
\[
(10),\ (9,1),\ (8,2),\ (8,1,1),\ (7,3),\ (6,4),\ (5,5)
\]
contain no complete curves.
\end{corollary}

\begin{proof}
Apply Proposition~\ref{prop:numerics} and Theorem~\ref{thm:strict}.
\end{proof}

The conclusions hold componentwise, including for the spin and hyperelliptic components classified in~\cite{KZ}.

As a complete subvariety of any positive dimension in a quasi projective variety would imply the existence of a complete curve, we obtain the following.  

\begin{corollary}\label{cor:subvarieties}
Every stratum in Corollaries~\ref{cor:lowgenus} and~\ref{cor:g6} contains no positive-dimensional complete subvariety.
\end{corollary}

\enlargethispage{6\baselineskip}

\end{document}